\documentclass[12pt]{amsart}

\usepackage{amsmath,amsthm,epsf}

\newtheorem{thm}{Theorem}[section]
\newtheorem{lem}[thm]{Lemma}
\newtheorem{cor}[thm]{Corollary}
\newtheorem{prop}[thm]{Proposition}
\newtheorem{rem}[thm]{Remark}

\newcommand{\bbr}{\begin{rem}\em} 
\newcommand{\eer}{\end{rem}}
\newcommand{\bpr}{\begin{proof}} 
\newcommand{\epr}{\end{proof}}

\def\Spin{\operatorname{Spin}}

\def\tb{\operatorname{tb}}

\def\Z{\hbox{$\mathbb Z$} }

\def\R{\hbox{$\mathbb R$} }

\def\tb{\operatorname{tb}}

\def\dfn#1{{\em #1}}

\def\tm#1{Theorem~\ref{#1}}

\def\Lm#1{Lemma~\ref{#1}}
\def\pr#1{Proposition~\ref{#1}}

\def\fig#1{Figure~\ref{#1}}

\begin{document}

\title{Tight Contact Structures on Lens Spaces}

\author{John B. Etnyre}
\address{Stanford University, Stanford, CA 94305}
\email{etnyre@math.stanford.edu}
\urladdr{http://math.stanford.edu/\char126 etnyre}


\begin{abstract}
	In this paper we develop a method for studying tight contact 
	structures on lens spaces.  We then derive uniqueness 
	and non-existence statements for tight contact structures with 
	certain (half) Euler classes on lens spaces.
\end{abstract}

\maketitle


\section{Introduction}

Contact geometry has recently come to the foreground of low 
dimensional topology.  Not only have there been striking advances in 
the understanding of contact structures on 3-manifolds 
\cite{el:twenty, gi, h}, but there has been significant interplay with 
knot theory \cite{r}, symplectic geometry \cite{lm, e:convex}, fluid dynamics 
\cite{eg:fluid}, foliation theory \cite{et}, and Seiberg-Witten theory 
\cite{km}.  In 1971 Martinet \cite{m} showed how to construct a 
contact structure on any 3-manifold. Several decades later 
it became clear that contact structures fell into two distinct classes:
tight and overtwisted (see Section~\ref{basic} for definitions).  It 
is the tight contact structures that carry significant geometric 
information; where as,
Eliashberg \cite{el:overtwisted} has shown that the understanding of overtwisted contact 
structures reduces to a ``simple'' algebraic question.  
Unfortunately, Martinet's theorem does not, in general, produce tight 
contact structures.  
The only general method for constructing
tight structures is by Stein fillings  \cite{g:stein, el:stein} or 
perturbing taut foliations \cite{et}.  These techniques, however, 
leave the general existence question open. Even less is known 
concerning the
uniqueness/classification question; it has only been answered on $S^3, S^{1}\times 
S^{2}$ and $\R P^{3}$ \cite{el:twenty} and $T^3$ (and most $T^{2}$ 
bundles over $S^{1}$)
\cite{gi, kanda}.

The purpose of this paper is to introduce some techniques for 
understanding tight contact structures.  We apply them to the simplest 
class of 
3-manifolds: lens spaces.  Recall lens spaces are 3-manifolds that
can be written as the union of two solid tori, or in other words, lens 
spaces are Heegaard genus one manifolds. On these manifolds we are able to derive 
some general uniqueness and non-existence statements in terms of the 
homotopy type of the contact structure. In particular, in 
Section~\ref{main-section} we show:

\begin{proof}[Theorem \ref{short}]
	{\em On any lens space $L(p,q)$ there is at least one class in 
    $H^{2}(L(p,q))$ realized by a unique tight contact structure and 
    at least one class that cannot be realized by a tight contact structure.}
\renewcommand{\qed}{}
\end{proof}
\renewcommand{\qed}{\qedsymbol}

It has been known for some time that on any 3-manifold there are only 
finitely many elements in second cohomology that can be realized by 
tight contact structures \cite{el:fill}.  This gives no 
restrictions for lens spaces since the second cohomology of a lens 
spaces is a finite 
group.  More recently Kronheimer and Mrowka \cite{km} have shown that only 
finitely many homotopy types of plane fields can be realized by 
semi-fillable 
contact structures.  Since any semi-fillable structure is tight and all currently known
tight structures are semi-fillable,
one is tempted to conjecture that this is also 
true for tight contact structures. This would restrict 
tight contact structures on lens spaces; it would not, however, say that there 
are a finite number of 
them. The techniques in this paper show:

\begin{proof}[Theorem \ref{finite}]
   {\em There are only finitely many tight contact structures (up to 
    isotopy) on any lens space.}
\renewcommand{\qed}{}
\end{proof}
\renewcommand{\qed}{\qedsymbol}

Moreover, on some lens spaces we can give a complete classification 
of contact structures.

\begin{proof}[Theorem \ref{complete}]
  {\em Classified up to isotopy:
  \begin{enumerate}
    \item If $p=0$ then $L(p,q)=S^{1}\times S^{2}$ and there is a 
    unique tight contact structure.
    \item If $p=1$ then $L(p,q)=S^{3}$ and there is a unique tight 
    contact structure.
    \item If $p=2$ then $L(p,q)=\R P^{3}$ and there is a unique tight 
    contact structure.
    \item On $L(3,1)$ there are exactly two tight contact structures 
    (one for each non zero element in $H^{2}(L(3,1);\Z)$).
    \item On $L(3,2)$ there is exactly one tight contact structure 
    (realizing the zero class in $H^{2}(L(3,2))$).
  \end{enumerate}}
\renewcommand{\qed}{}
\end{proof}
\renewcommand{\qed}{\qedsymbol}

In future work we plan to push the analysis of tight contact 
structures on lens spaces further
\footnote{
Added in proof: E.~Giroux and K.~Honda (independently) have recently announced a complete classification
of tight contact structures on lens spaces. Specifically they are all obtained from Stein fillings
and determined by their half-Euler class.
}  
as well as apply these techniques to 
3-manifolds with higher Heegaard genus.

\section{Contact Structures In Three Dimensions}\label{basic}

In this section we briefly recall some facts concerning contact 
geometry in dimension three.  For more details see \cite{a:etall} or 
\cite{el:twenty}.

Let $M$ be an oriented 3-manifold.
A \dfn{contact structure} $\xi$ on $M$ is a totally non-integrable
2-plane field in $TM.$  
We will only consider transversely orientable contact structures (this 
is not a serious restriction).  This allows us to globally define the 
plane field $\xi$ as the kernel of a 1-form $\alpha.$  Using the 
1-form $\alpha,$ the Frobenius Theorem allows us to express the
non-integrability of $\xi$ as $\alpha\wedge 
d\alpha\not=0.$ Thus 
$\alpha\wedge d\alpha$ is a volume form on $M.$  Other 1-forms we 
could use to define $\xi$ will give us different volume forms 
but they will induce the same orientation on $M$ (since 
any other 1-form $\alpha'$ that defines $\xi$ must be of the form
$f\alpha$ where $f:M\to \R$ is a non-zero function).  Thus $\xi$ 
orients $M$ and we only consider contact structures $\xi$ whose 
induced orientation agrees with the orientation on $M.$ (A similar analysis 
could be done when the orientations disagree.)

Contact geometry presents interesting and
difficult global problems; however, Darboux's Theorem tells us that
all contact structures are {\em locally} contactomorphic. Two contact
structures are \dfn{contactomorphic} if there is a diffeomorphism of the
underlying manifolds that sends one of the plane fields to the other.  
Furthermore, Gray's Theorem tells us that a continuously varying 
family of contact structures are related by a continuously varying 
family of contactomorphisms.
We have similar results near surfaces in $M.$
If $\Sigma$ is a surface in a contact manifold $(M,\xi)$ then generically
$T_p\Sigma\cap \xi_p$ will be a line in $T_p\Sigma.$  
Since a line field is always integrable $T_p\Sigma\cap \xi_p$
defines a natural
singular foliation $\Sigma_\xi$ associated to $\xi$ called the 
\dfn{characteristic foliation}.  As with Darboux's Theorem determining
a contact structure in the neighborhood of a point, one can show that 
$\Sigma_\xi$ determines the germ of $\xi$ along $\Sigma.$
The singular points of $\Sigma_{\xi}$ occur where $\xi_{p}=T_{p}\Sigma.$  
Generically, singular points of $\Sigma_{\xi}$ are either elliptic 
(if the local index is 1) or hyperbolic (if the local index is $-1$).
Notice that if $\Sigma$ is oriented then we can assign signs to the 
singular points of $\Sigma_{\xi}.$  A singular point $p$ is called 
\dfn{positive} if $\xi_{p}$ and $T_{p}\Sigma$ have the same 
orientation, otherwise $p$ is called \dfn{negative}.  Moreover, the 
orientations on $\Sigma$ and $\xi$ ``orient'' the characteristic 
foliation, so we can think of the singular foliation as a flow.
A very 
useful modification of $\Sigma_{\xi}$ is given by the Elimination 
Lemma (proved in various forms by
Giroux, Eliashberg and Fuchs, see \cite{el:knots}).

\begin{lem}[Elimination Lemma]\label{lem:elimination}
	Let $\Sigma$ be a surface in a contact $3$-manifold $(M,\xi)$.  Assume
	that $p$ is an elliptic and $q$ is a hyperbolic singular point in 
	$\Sigma_\xi$, they both have the same sign and there is a leaf $\gamma$
	in the characteristic foliation $\Sigma_\xi$ that connects $p$ to $q$.
	Then there is a $C^0$-small isotopy $\phi:\Sigma\times[0,1]\to M$ such that
	$\phi_0$ is the inclusion map, $\phi_t$ is fixed on $\gamma$ and outside any 
	(arbitrarily small) pre-assigned
	neighborhood $U$ of $\gamma$ and $\Sigma'=\phi_1(\Sigma)$ has no singularities
	inside $U$.
\end{lem}

There has recently emerged a fundamental dichotomy in three dimensional 
contact geometry.
A contact structure $\xi$ is called \dfn{overtwisted} if there exists an
embedded disk $D$ in $M$ whose characteristic foliation $D_\xi$ contains 
a limit cycle.  If $\xi$ is not overtwisted then it is called
\dfn{tight}.   Eliashberg \cite{el:overtwisted} has completely
classified overtwisted contact structures on  closed 3-manifolds:  
classifying overtwisted contact structures up to isotopy is 
equivalent to classifying plane fields up to homotopy (which has a 
purely algebraic solution).
As discussed in the introduction, much less is known about tight 
contact structures.  One of the main results about tight contact 
structures, on 
which all the results in this paper are based, is the following 
theorem of Eliashberg.

\begin{thm}[Eliashberg \cite{el:twenty}, 1992]\label{thm:contact.ball}
	Two tight contact structures on the ball $B^3$ which induce the same
	characteristic foliations on $\partial B^3$ are isotopic relative
	to $\partial B^3$.
\end{thm}

A closed curve $\gamma:S^{1}\to M$ in a contact manifold $(M,\xi)$ is called 
\dfn{transverse} if $\gamma'(t)$ is transverse to $\xi_{\gamma(t)}$ 
for all $t\in S^{1}.$  
It can be shown that any curve can be made 
transverse by a $C^{0}$ small isotopy.
Notice a transverse curve can be 
\dfn{positive} or \dfn{negative} according as $\gamma'(t)$ agrees with 
the co-orientation of $\xi$ or not. We will restrict our attention to 
positive transverse knots (thus in this paper ``transverse'' means 
``positive transverse''). 
Given a transverse knot $\gamma$ in $(M,\xi)$ that bounds a surface $\Sigma$
we define the \dfn{ self linking number}, $l(\gamma)$, of $\gamma$ as follows: 
take a non-vanishing vector field $v$ in $\xi\vert_{\gamma}$ that 
extends to a non-vanishing vector field in $\xi\vert_{\Sigma}$ and let $\gamma'$
be $\gamma$ slightly pushed along $v$. Define
$$l(\gamma,\Sigma)=I(\gamma',\Sigma),$$
where $I(\cdot,\cdot)$ is the oriented intersection number.
There is a nice relationship between $l(\gamma,\Sigma)$ and the 
singularities of the characteristic foliation of $\Sigma.$  Let 
$d_{\pm}=e_{\pm}-h_{\pm}$ where $e_{\pm}$ and $h_{\pm}$ are the 
number of $\pm$ elliptic and hyperbolic points in the characteristic 
foliation $\Sigma_{\xi}$ of $\Sigma,$ respectively.  In \cite{be} it 
was shown that
$$l=d_{-}-d_{+}.$$
When $\xi$ is a {\em tight} contact structure and $\Sigma$ is a 
{\em disk},
Eliashberg \cite{el:twenty} has shown, using the elimination lemma,  how to 
eliminate all the positive hyperbolic and negative elliptic points 
from $\Sigma_{\xi}$ (cf.\ \cite{gi:convev}). Thus in a 
tight contact structure $l(\gamma,\Sigma)$ is always negative.

\section{Singular Foliations on the 2-Skeleton}

The lens space $L(p,q)$ is, by definition, the union of two solid tori
$V_{0}$ and $V_{1}$ glued together by a map $\phi:\partial V_{1}\to \partial V_{0}$ which in standard 
coordinates on the torus is given by the matrix
$$\begin{pmatrix} -q&p'\\ p&r\end{pmatrix},$$
where $p'$ and $r$ satisfy $-rq-pp'=-1$.
(A standard basis is given by 
$\mu$, the boundary of a meridianal disk, and $\lambda$, a longitude 
for $T^2$ given by the product structure on $V_i$ and oriented so that 
$\mu\cap\lambda=1$.) We will also use a CW decomposition of $L(p,q)$ obtained
from $V_{0}$ and $V_{1}$ as follows: let $C$ be the core curve in $V_{0},$
$x$ a point on $C,$ $D$ a disk in $L(p,q)$ that intersects $V_{1}$ in  a
meridianal disk and whose boundary wraps $p$ times around $C$, and 
$B=L(p,q)\setminus D$ a 3-ball. Now $L(p,q)$ can be written as
$$\{ x\}\cup C \cup D\cup B.$$
We call $D$ the \dfn{generalized projective plane} in $L(p,q).$ Note 
that if $p=2$ so that $L(p,q)$ is $\R P^{3}$ then $D$
is a copy of $\R P^{2}\subset L(p,q).$  Given a 
contact structure $\xi$ on $L(p,q)$ we would like to understand the
characteristic foliation, $D_{\xi},$ on $D.$  To this end we begin by 
isotoping $C$ so that it is transverse to $\xi.$  
Throughout this paper we will take $V_0$ to be a small standard neighborhood of the transverse
curve $C.$
We now have 
the following proposition which says that $\xi$ on all of $L(p,q)$  is 
determined by $D_{\xi}.$

\begin{prop}\label{foliation-determines}
	Let $L_0$ and $L_1$ be two copies of $L(p,q)$. Let $\xi_i$ be a tight oriented
	contact structure on $L_i$ and $D_i$ the generalized projective plane
	in $L_i$, $i=0,1$.  Assume that the 1-skeleton $C_i$ of $D_i$ is 
	transverse to $\xi_i$. If a diffeomorphism $f:L_0\to L_1$ may be isotoped so that
	it takes $(D_0)_{\xi_0}$ to $(D_1)_{\xi_1}$ then it may be isotoped 
	into a contactomorphism.
\end{prop}

\bpr
A slight modification of the standard proof that the characteristic 
foliation on a surface determines the contact structure on a 
neighborhood of the surface can be used to isotope $f$ into a 
contactomorphism in a neighborhood of $D_{0}.$  We can then use 
Eliashberg's characterization of tight contact structures on the
3-ball, Theorem~\ref{thm:contact.ball}, to isotope $f$ into a contactomorphism 
on all of $L_{0},$ since $L_{0}\setminus D_{0}$ is a 3-ball and $f$ is already a 
contactomorphism in a neighborhood of its boundary.
\epr

In the remainder of this section we derive a standard form for the 
characteristic foliation on the two skeleton $D\subset L(p,q)$ depending 
only on the homotopy class of $\xi.$

\subsection{The Euler class and the Singular Foliation}

We have already arranged that $C$ is transverse to $\xi;$ thus, if $V_{0}$ is a sufficiently 
small tubular neighborhood of $C,$ the curve 
$\gamma=\partial V_{0}\cap D$ will also be a transverse curve.  It 
will be helpful to keep in mind that $\gamma$ is homotopic to $pC.$ We 
will interchangeably think of $D$ as an embedded disk with boundary 
$\gamma$ and the generalized projective plane in $L(p,q)$ (note this 
should cause no confusion since one uniquely determines the other).

As discussed in Section~\ref{basic} the self-linking number 
$l=l(\gamma,D)$ of $\gamma$ is related to the 
singularities of $D_{\xi}.$  
We can also relate $l=l(\gamma,D)$ to the homotopy class of $\xi.$  To 
state this relation we must begin by recalling the definition of a 
$\Z_{2}$-refinement of the Euler class of $\xi.$  
The Euler class $e=e(\xi)$ is
an even element of $H^{2}(L(p,q);\Z)=\Z_{p}.$  If $p$ is odd then then
there is a unique element $c$ in $H^{2}(L(p,q);\Z)$ such that $2c=e.$  
However, if $p$ is even then there are precisely two elements in 
$H^{2}(L(p,q);\Z)$ that can be thought of as half of $e.$  We cannot naturally associate
one of these elements to $\xi,$ however once a spin structure is fixed on $L(p,q)$ we
can. Thus the 
$\Z_{2}$-refinement $\Gamma(\xi)$ of $e(\xi)$ is actually a map
$$\Gamma(\xi):\Spin(L(p,q))\to G,$$
where $\Spin(L(p,q))$ is the group of spin structures on $L(p,q)$ and 
$G=\{x\in H_1(L(p,q);\Z): 2x=PD(e(\xi))\}$ ($PD$ means Poincar\'e dual). 
For a general discussion 
of this invariant see \cite{g:stein} and for an
exposition of spin structures see \cite{gs}.

Recall, if $p$ is odd then $L(p,q)$ has a unique spin structure and
if $p$ is even there are precisely two and they may be distinguished 
as follows:  given a spin structure on $L(p,q)$ there is a naturally induced spin 
structure on $L(p,q)\times [0,1]$ and one may ask if this spin 
structure extends over a 4-dimensional 2-handle added to the curve $C$ 
in $L(p,q)\times \{1\}.$  One of the spin structures on $L(p,q)$ will 
extend over a 2-handle added with framing 0 and one will not.
Note when $p$ is odd $\Gamma(\xi)$ is determined by the element in 
$H^{2}(L(p,q);\Z)$ that is half of $e(\xi)$ thus we refer to 
$\Gamma(\xi)$ as the \dfn{half-Euler class 
of $\xi$}.
In general, $\Gamma(\xi)$ clearly refines 
the Euler class since $2\Gamma(\xi)(s)=PD(e(\xi))$ for
any $s\in \Spin(L(p,q))$. 

The map $\Gamma(\xi)$ is one-to-one and for 
computational convenience we shall actually define it as a map from 
$G$ to $\Spin(L(p,q)).$  For an alternate definition and its relation to 
the one we give here, in particular its well-definedness, see 
\cite{g:stein, e:dis}.
To define $\Gamma(\xi)$ 
let $v$ be a vector field in $\xi$ with zero locus (counted with multiplicity) $2c$
where $c$ is a smooth curve in $L(p,q)$.  Note that the homology class of $c$ is in $G.$ 
The vector field $v$ gives a trivialization
of $\xi$ on $L(p,q)\setminus c$ and hence a trivialization of $TL(p,q)$ on 
$L(pq)\setminus c$.
This trivialization induces a spin structure on $L(p,q)\setminus c$, and finally,
since $v$ vanishes to order two on $c$ we can extend this spin structure over
$c$ to obtain a spin structure $s'$ on all of $L(p,q)$.  The map $\Gamma(\xi)$ will associate $c$ and $s'.$
Specifically we define
$$\Gamma(\xi)(s')=[c].$$    We can now state the relation between $l(\gamma)$ 
and $\Gamma(\xi).$

\begin{thm}\label{l-formula}
	Let $s'$ be any spin structure on $L(p,q).$ If $p>0$ is even, 
	let $s$ be the spin structure on $L(p,q)$ that does not
	extend over a 2-handle attached to $C$ with framing 0.
    	Then 
	\begin{equation}\label{better.l-formula}
		\Gamma(\xi)(s')\cdot D\equiv\frac{1}{2}(-l(\gamma)+q+p\Delta(s,s'))\mod p,
	\end{equation}
	where $\Delta(s,s')=0$ if $s=s'$ and $1$ otherwise.
	If $p$ is odd this formula also holds if we always take $\Delta (s,s')$ to be 1 when $q$ is even
	and 0 when $q$ is odd.
\end{thm}

\bbr
   When $p$ is odd $L(p,q)$ has a unique spin structure which, when $q$ 
   is even, is not the structure
   $s$ described in the theorem. This explains then strange definition of $\Delta (s,s')$ 
   for odd $p$.
\eer

\bbr\label{oddp}
   Note that in terms of the Euler class we get
   	\begin{equation}\label{lformula}
		e(\xi)(D)\equiv(-l+q)\mod p.
	\end{equation}
   This formula is easier to prove and would suffice for $p$ odd. 
\eer

\bpr
   We begin by constructing a vector field $w$ in $\xi\vert_{D}$
   and use $w$ to compute the Euler class of $\xi.$ 
   Then we modify $w$ to a vector field whose zeros all have 
   multiplicity two, thus allowing us to compute $\Gamma.$ On $D\cap 
   V_{1}$ we let $w$ be a vector field directing the characteristic 
   foliation $D_{\xi}.$  We cannot use the characteristic foliation 
   to define $w$ on $D\cap V_{0}$ since it will not be well defined 
   along $C.$  But it is not hard to find a vector field in $\xi\vert_{D\cap 
   V_{0}}$ that agrees with $w$ on $\xi\vert_{\partial (D\cap V_{0})}$ and has 
   exactly $q$ zeros in $V_{0}.$  (If $p$ is even this vector field will define a 
   spin structure in a neighborhood of $C$ that will not extend over 
   a 2-handle attached to $C$ with framing 0.) When one uses $w$ to compute the 
   Euler class of $\xi$ one gets $-l+q.$
   
   We now must coalesce the zeros of $w$ into zeros with multiplicity two.
   For simplicity assume that $p$ is even since otherwise $\Gamma$ is 
   determined by $e$ and we are already done.  
   We can also assume that we have canceled all the negative elliptic 
   points and positive hyperbolic points in $D_{\xi},$
   thus there will be $n+1$ elliptic points and $n$ hyperbolic points
   (see the remark following this proof).
   Writing down a model 
   for $D\cap V_{0}$ one can explicitly write a vector field with 
   $\frac{q-1}{2}$ zeros of multiplicity 2 and one zero of multiplicity 
   1.  Moreover this vector field will agree with $w$ on $C$ and 
   $\partial (D\cap V_{0}).$ Let $v$ be this vector field on $V_{0}.$
   The zero of multiplicity 1 in $V_{0}$ will be connected by a leaf 
   to an elliptic point on $D\cap V_{1}$ and the other elliptic and 
   hyperbolic points will all pair up along stable separatrices of 
   the hyperbolic points.  It is now not hard to explicitly write down a
   vector field in a neighborhood of the connecting leaves of these 
   pairs that agree with $w$ outside the neighborhood and has 
   precisely one zero of
   multiplicity 2 on each seperatrix.  This will allow us to 
   define $v$ on all of $D.$  One can now extend $v$ to all of 
   $L(p,q)$ and as mentioned above the spin structure it induces is 
   $s.$  Thus we have $\Gamma(\xi)(s)\cdot D\equiv\frac{1}{2}(-l+q)\mod 
   p.$  The formula for the other spin structure follows from general 
   properties of $\Gamma$ (see \cite{g:stein}).
\epr

\bbr\label{eh-formula}
Theorem~\ref{l-formula} tells us that
\begin{equation} \label{prac-l}
l(\gamma)=q-2(\Gamma(\xi)(s)\cdot D)+2np.
\end{equation}
Recall, when $\xi$ is tight $l$ must be negative,  limiting the 
possibilities for $n.$
Moreover we can cancel all the negative 
elliptic and 
positive hyperbolic singularities from the characteristic foliation 
$D_{\xi}.$ Thus we have
\begin{equation}
e_{+}+h_{-}=-q+2(\Gamma(\xi)(s)\cdot D)+ 2np.
\end{equation}
We also know that $e_{+}-h_{-}=1$ since we can use the characteristic 
foliation $D_{\xi}$ to compute the Euler characteristic of the disk 
$D.$  In the remainder of this section we will simplify the 
characteristic foliation more, eventually showing $0<e_{+}+h_{-}< 2p.$  
Notice that this will uniquely determine $e_{+}$ and $h_{-}.$ (Again, 
when $p$ is odd and $q$ is even one should use Equation~\eqref{lformula} 
as discussed in Remark~\ref{oddp}.)
\eer

\subsection{Making Stars}

The \dfn{graph of singularities} of $D_{\xi}$ will be defined to be 
the union of all singular points and stable separatrices (i.e. points 
that limit to a hyperbolic singularity in forward time) in $D_{\xi}.$
When talking of the graph of singularities we  always 
assume that $e_{-}=h_{+}=0.$
It is easy to see that the graph of singularities must be a tree, but 
one can say much more.

\begin{figure}[ht]
	{\centerline{\epsfbox{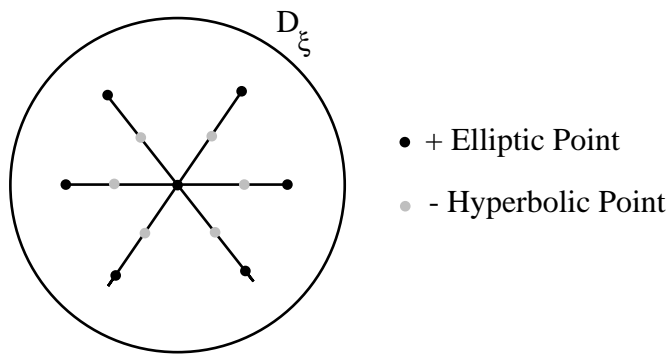}}}
	\caption{Singularities on $D.$}
	\label{fig:star}
\end{figure}

\begin{lem}\label{star}
	We can choose  $D$ so that the graph of singularities of $D_\xi$ 
	forms a star, i.e.\  there is one $(e_{+}-1)$-valent elliptic vertex, ($e_{+}-1$) univalent
	elliptic vertices and exactly one hyperbolic singularity in the interior of each edge
	(see \fig{fig:star}).
\end{lem}

\bpr
This is a special case of a lemma in \cite{ml}, though the proof there seems
to be incomplete. Some of the ideas below are also reminiscent of ones appearing in 
Fraser's thesis \cite{f}, though in a different setting.

\begin{figure}[ht]
		{\centerline{\epsfbox{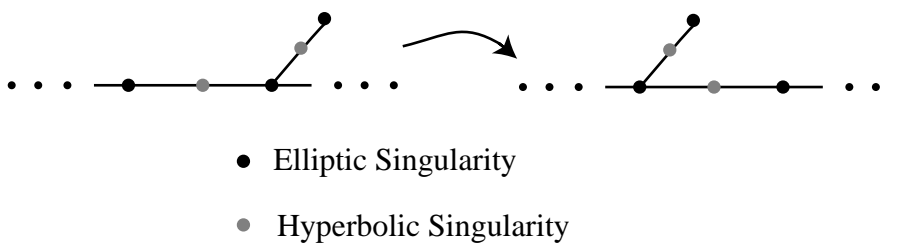}}}
		\caption{Change in Graph of Singularities.}
		\label{fig:sing.change}
	\end{figure}

We will show how to isotope $D$ to a disk $D'$ with  transverse boundary
in $\partial V_1,$ whose  graph of singularities relates to $D$'s as 
shown in Figure~\ref{fig:sing.change}.  Since the graph of singularities in 
$D_\xi$ must be a tree, a sequence of such moves will clearly yield the 
conclusion of the lemma. 
Assume that part of the graph of singularities in $D_\xi$ is as shown 
on the left hand side of 
Figure~\ref{fig:sing.change}.  Let $D_i$ be a subdisk of $D$ with transverse boundary containing 
the graph of singularities in $D_\xi.$
Let $U$ be a neighborhood, diffeomorphic to an open ball, of $D_i$ in $V_1$ and set $D_a=D\cap U.$ 
Now $(U, \xi\vert_U)$ is a tight contact structure on $\R^3$, so the classification of contact structures
on $\R^3$ \cite{el:knots} implies there is a contactomorphism
$$f:(U,\xi\vert_U)\to (\R^3,\xi'=\{dz+xdy=0\}).$$  
In the following paragraph we show that there is 
a compactly supported isotopy of $f(D_a)$ to $D_a'$ so that $(D'_a)_{\xi'}$ is related to $D_\xi$ as shown in 
Figure~\ref{fig:sing.change}. Then our desired disk is $D'= (D\setminus D_a)\cup f^{-1}(D'_a)$ (note the
two pieces fit together since the above isotopy was compactly supported and the open disk $f(D_a)$ is
properly embedded in $\R^3$).

\begin{figure}[ht]
	{\centerline{\epsfbox{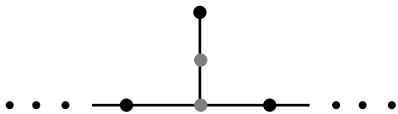}}}
	\caption{Characteristic Foliation on $\Delta$.}
	\label{fig:graph.model}
\end{figure}

We now need to prove our claim concerning the compactly supported isotopy of $f(D_a).$ 
To this end let $\Delta$ be a disk in the standard contact structure on $\R^3$ whose characteristic foliation
is related to $D_\xi$ as indicated in Figure~\ref{fig:graph.model} (creating such a disk is an easy 
exercise).
Note that the characteristic foliation on $\Delta$ is unstable. Specifically, if we take a point $p$ on 
the seperatrix connecting the two hyperbolic points and push it up, respectively down, $\Delta_\xi$ will look like
the foliation indicated on the right side, respectively left side, of Figure~\ref{fig:sing.change}.
For later use we describe this isotopy in a standard model: let $p$ be a point on the seperatrix, 
$\alpha$ a segment of the seperatrix containing $p,$ 
and in $(\R^3,\xi')$ let $p'$ be the point $(1,0,0)$ and $\alpha'$ the arc $\{(t,0,0): \frac{1}{4}\leq t\leq 2\}.$
Now there is a neighborhood $O$ of $\alpha$ diffeomorphic
to an open three ball and contactomorphic to a neighborhood 
$O'$ of $\alpha'$ that does not intersect the $yz$-plane in $\R^3.$ 
Denote this contactomorphism by $g.$ Moreover, if we take $O$ and $O'$ sufficiently small we can
assume that $g$ takes $O\cap \Delta$ to the $O'\cap (xy$-plane$),$ $\alpha$ to $\alpha'$ and $p$ to $p'.$ 
One may now explicitly see that
slightly pushing $p'$ up or down (in the $z$-direction) 
will change the foliation as claimed. If $O'$ is an $\epsilon$ neighborhood of
$\alpha'$ (which we may assume by shrinking $O$ and $O'$) 
then choose a $\delta<<\epsilon$ and push $p'$ up slightly by an isotopy supported in the
$\delta$-ball around $p'.$  Let $\Delta'$ be the image of $\Delta$ under the corresponding isotopy.
Below we will find a compactly supported contactomorphism $h:\R^3\to \R^3$ taking $f(D_i)$ to $\Delta'$ and
a neighborhood $N$ of $f(D_i)$ to a neighborhood $N'$ of $\Delta'.$ 
Assuming this for the moment we finish the proof of the lemma. 
By construction, $h(f(D_a))$ has the same characteristic foliation as $D.$ To obtain the
desired foliation  
we can isotope $h(f(D_a))$ in $O$ to a disk we call $D'_a$ so its characteristic
foliation is as shown on the right of Figure~\ref{fig:sing.change} (this isotopy corresponds to pushing
$p'$ down a little in $O'$). 

The problem now is that $D'_a$ may not
be embedded if  
$$K=h(f(D_a\setminus D_i))\cap O\not=\emptyset.$$ 
Consider $K'=g(K)$ in our neighborhood
$O'$ of $p'$ and by possibly shrinking $N$ and $N'$ we can assume that $g(N')$ is a $\delta'$ neighborhood
of $g(\Delta'\cap O)$ where $\delta'<< \delta.$  Finally take $m$ to be the maximum $z$ coordinate
in $g(\Delta'\cap O)$ (though not essential we can assume that $g(\Delta'\cap O)$ differs from 
the $xy$-plane by a small symmetric bump with a unique maximum $m$). Note that if $\delta'\geq m$ then 
$D'_a$ will be embedded because the isotopy from $h(f(D_a))$ to $D'_a$ does not have to leave
the neighborhood $N'.$ 
So we assume that this
is not the case. We can now choose a function $\psi:\R\to\R$
such that 
\begin{enumerate}
\item{} $\psi(z)=z$ outside $[-\delta', m],$ 
\item{} $\psi(z)<0$ for $z<m-\delta',$ 
\item{} $\frac{1}{2}\leq \psi'(z)\leq 1$ on $[-\delta',m-\delta']$  and 
\item{} $\psi$ is strictly increasing. 
\end{enumerate}
With $\psi$ in hand we define a contactomorphism $\phi:\R^3\to \R^3$ by $\phi(x,y,z)=(f'(z)x, y f(z)).$
So $\phi(K')$ clearly has the same characteristic foliation as $K'$ (in particular it is non singular) 
and agrees with $K'$ near the boundary of $O'.$ 
Using properties 2.\ and 3.\ of $\psi$ we see that $\phi(K')$ lies in $O'$ and below the $xy$-plane.
Thus if we replace $K'$ by $\phi(K')$ (and of course make
the corresponding alteration to $D'_a$) then we have eliminated all the possible self intersections of $D'_a.$

To finish the proof we now need to show that
there is a compactly supported contactomorphism $h:\R^3\to \R^3$ taking a neighborhood $N$ of $f(D_i)$ to 
a neighborhood $N'$ of $\Delta'.$ 
Recall that as contact manifolds $\R^3=S^3\setminus\{q\},$ where $S^3$ has
its standard tight contact structure and $q$ is any point in $S^3.$ 
So we actually construct a contactomorphism 
of $S^3$ that fixes a neighborhood
of $q.$ Since $f(D_i)$ and $\Delta'$ have identical characteristic foliations we can find a diffeomorphism $h$ of $S^3$
fixing a neighborhood of $q,$ taking $f(D_i)$ to $\Delta'$ and restricting to a contactomorphism on the 
neighborhoods $N$ and 
$U_0$  of $f(D_i)$ and $q,$ respectively. We may now isotop $g$ into a contactomorphism on all of $S^3$ (this easily 
follows from
the proof of Eliashberg's classification of tight contact structures on $B^3$ \cite{el:twenty} since $g$ gives
a map from the 3-ball $S^3\setminus N$
to the 3-ball $S^3\setminus g(N)$ preserving the characteristic foliation on their boundaries and 
restricting to a contactomorphism
on the Darboux ball $U_0$). 
\epr

\bbr\label{setup}
    We  now set up some notation that will be used through the 
    rest of the paper. Assume that we have already arranged that the 
    graph of singularities of $D_{\xi}$ is a star. If $n+1$ is the number 
    of elliptic points in 
    $D_{\xi}$ then there are $n$ hyperbolic points which we 
    denote $h_{1},\ldots,h_{n}$ numbered anti-clockwise.
    We also label the 
    $n$-valent elliptic point $e_{0}$ and the other elliptic points $e_{i}$ 
	according to the hyperbolic point with which they share a flow line.
   	A point $x$ on the 1-skeleton $C$ will break the boundary of $D$ 
	into $p$ arcs, $B_1, B_2,\ldots, B_{p}$ (also numbered 
	anti-clockwise).  Notice that as we traverse the boundary of $D$ 
	anti-clockwise we will encounter both the end points of the unstable 
	separatrices (the ones {\it not} shown in Figure~\ref{fig:star}) 
	leaving $h_{1}$ then both leaving $h_{2}$ continuing in this fashion 
	until we reach the end points coming from $h_{n}$.  We 
	label the end points of the unstable separatrices leaving $h_{i}$
	as $h_{i}^{a}$ and $h_{i}^{c}$ so that $h_{i}^{a}$ is anti-clockwise 
	of $h_{i}^{c}.$ 
	Thus around the boundary of $D$ we see the points 
	$h_{1}^{c},h_{1}^{a},h_{2}^{c},\ldots, h_{n}^{c},h_{n}^{a}$  broken into 
	$p$ sets by the intervals $B_{i}.$
\eer

\subsection{Simplifying Stars}

So far we have arranged that the generalized projective 
plane $D$ has the following properties:
\begin{itemize}
\item the one skeleton of $D$ is transverse to $\xi,$
\item there are no negative elliptic and no positive hyperbolic 
singularities in $D_{\xi},$ and
\item the graph of singularities in $D_{\xi}$ forms a star.
\end{itemize}
We can now get control of the number of branches in the graph of 
singularities.

\begin{lem}\label{simplify}
    If $e_{+}>p$ then we can isotope $D$ to $D'$ so that $D'$ 
    is a generalized projective plane in $L(p,q)$ enjoying the above 
    listed properties for $D$ and $e_{+}(D')=e_{+}(D)-p.$
\end{lem}

\bpr
	We know the singularities of $D_\xi$ form
	a star.  Since we are assuming that $e_{+}>p$ we know there are at least
	$p$ edges in the star, hence there are $n\geq p$ hyperbolic points 
	(one for each edge). Now using the notation of Remark~\ref{setup} it 
	is clear that if $n>p$  then at least one of the hyperbolic points,
	$h_{i}$ say, has both its unstable separatrices exiting $D$ through, 
	say, $B_{j}.$
        This is also true when $n=p.$ To 
	see this choose the point $x$ on $C$ so that 
	$h_{1}^{c}$ is the closest point (in the anti-clockwise direction) to 
	$x$ and lying in $B_{1}.$
	Now if $h_{1}$ is not the point we are looking for then $h_{1}^{a}$ 
	must be leaving $D$ through $B_{2}$ (or an interval further 
	anti-clockwise).  If we continue in this fashion and none of the points
	$h_{i}$ for $i<p$ have their unstable separatrices leaving on the 
	same $B_{j}$ then $h_{p}^{a}$ and $h_{p}^{c}$ must both lie in 
	$B_{p},$ thus proving our claim.
	
\begin{figure}[ht]
		{\centerline{\epsfbox{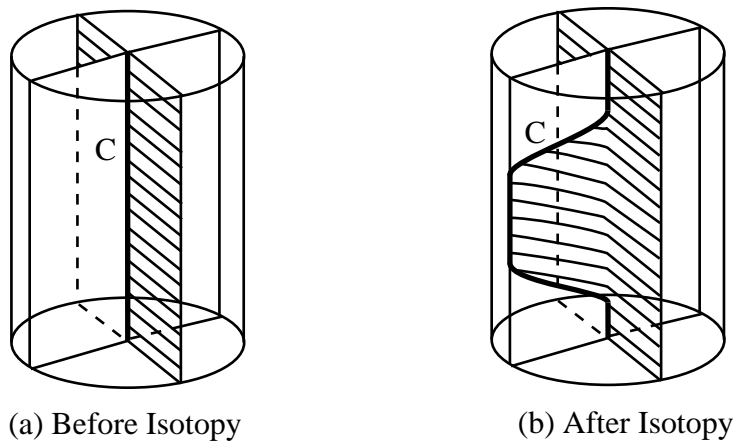}}}
		\caption{Model for $V_0.$}
		\label{fig:loc.model}
	\end{figure}

\begin{figure}[ht]
		{\epsfxsize=5truein\centerline{\epsfbox{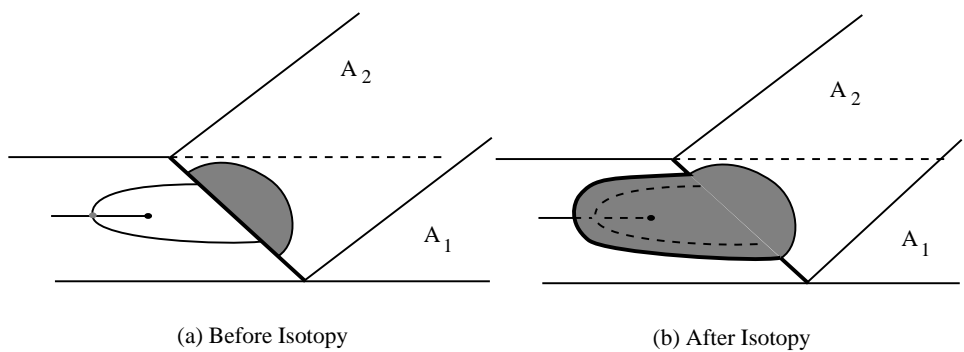}}}
		\caption{$A_1$ and $A_2$ near $C\cap V_0.$}
		\label{fig:fan.iso}
	\end{figure}

	Having found this hyperbolic point $h_{i}$ with both unstable 
	separatrices leaving $D$ through $B_{j}$ we now describe an isotopy 
	of $D$ which will decrease $e_{+}$ by $p.$  Note that the unstable 
	separatrices of $h_{i}$ separate a disk $\Delta$ from $D$ which contains
	exactly one elliptic point $e_{i}$ and has part of its boundary on $B_{j}$ 
	and the other part is made from the unstable separatrices of $h_{i}.$  
	We use this disk $\Delta$ to guide our isotopy.
	Essentially we push (in an 
	arbitrarily small neighborhood of $\Delta$) the part of $C$ that 
	intersects $\Delta$ across the 
	unstable separatrices of $h_{i}.$    More precisely we can write 
	down an exact model of our situation in $V_{0}$ and then isotope the 
	interval 
	$I=C\cap \Delta$ to the boundary of $V_{0}$ along $\Delta$ (see 
	Figure~\ref{fig:loc.model}).
	Next we isotope $I$ across the
	unstable separatrices of $h$.  One can write down a 
	precise local model of $\Delta$ in which 
	to do the isotopy.
	Note this is not an isotopy through transverse
	knots, as the one above, but at the end of the isotopy $C$ will again
	be a transverse knot. To see what happens to the rest of $D$ we 
	consider the case when $p=3$,  the case for larger $p$ being analogous. 
	Near $I$ we have $\Delta$ on one side of $I$ and the other two branches 
	of $D\cap V_{0},$ which we label $A_{1}$ and $A_{2},$ 
	fanning out behind $C$ (see \fig{fig:fan.iso} (a)).  
	We can assume that $A_{1}$ meets $I$ so that it and $\Delta$ form
	a smooth surface.  When we push $I$
	across the unstable separatrices of $h_{i}$ we will transfer $h_{i}$ and $e_{i}$ 
	from $\Delta$
	to $A_{1}.$
	But the orientation $\Delta$  inherits 
	after the isotopy, i.e.\ as a subset of $A_{1},$ is opposite the   
	orientation it originally inherited 
	from $D$ (to see this consider the situation when $p=2$ and we are 
	dealing with a projective plane).
	We will also have to drag part
	of $A_{2}$ along with us through the isotopy 
	(the gray part in \fig{fig:fan.iso}), 
	but notice that we can drag it so that
	it is arbitrarily close to $A_1$ (see \fig{fig:fan.iso} (b)).  Thus 
	the characteristic foliation on the part of $A_2$ we dragged along is topologically
	equivalent to the foliation just transferred onto $A_1,$ (since the foliation
	is structurally stable).  
	In particular, this means that we have added an elliptic
	and a hyperbolic singularity to $A_2$.  After isotoping $C$ we have a new
	transverse curve $C'$ and a new generalized projective plane $D'$.  Using
	$C'$ we also get a new Heegaard decomposition $L=V_0'\cup V_1'$ where
	$V_0'$ is a small tubular neighborhood of $C'$.
	
\begin{figure}[ht]
		{\centerline{\epsfbox{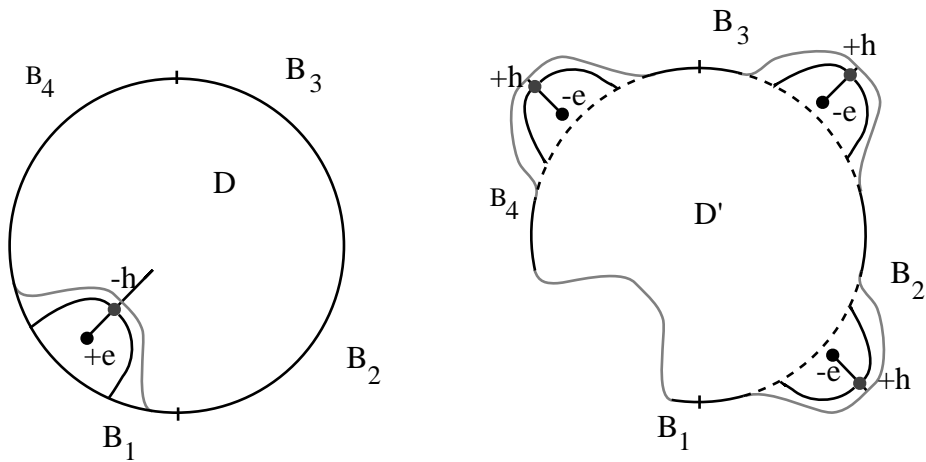}}}
		\caption{$D$ and $D'$ when $p=4$.}
		\label{fig:newd}
	\end{figure}

	We now claim that after canceling all the newly created negative 
	elliptic and positive hyperbolic points on $D'$ then the number of 
	positive elliptic points, $e_{+}',$ is $e_{+}-p.$
	To see this note
	$D'$ will essentially
	look like $D$ with $\Delta$ removed in one place and $p-1$ copies of
	it glued on along subarcs of $B_2,\ldots,B_{p}$ (see \fig{fig:newd}).
	The orientations on the copies
	of $\Delta$ in $D'$ will be opposite that of $\Delta$ in $D$.  Thus
	$D'_\xi$ has one less positive elliptic point and one less negative hyperbolic
	point than $D_\xi$ but has $p-1$ more negative elliptic and positive 
	hyperbolic points.  Thus when we cancel the negative elliptic and positive
	hyperbolic points from $D'_\xi$ we will have $p$ fewer positive elliptic
	and negative hyperbolic point than $D_\xi$ had. 
\epr 

\bbr\label{big-point}
  Notice that the proof shows that, in a tight contact structure, 
  if one is already in a minimal 
  configuration, i.e.\ $e_{+}\leq p,$ then the unstable separatrices 
  emanating from one hyperbolic point cannot both leave $D$ along the 
  same arc $B_{i}$ no mater which point $x$ is used to form the arcs 
  $B_{i}.$  This observation will be crucial in what follows. 
\eer

\section{Contact Structures on Lens Spaces}

\subsection{The Existence of Tight Contact Structures}

The main technique for generating tight contact structures on 
3-manifolds is to realize the 3-manifold at the ``boundary'' of a 
Stein manifold \cite{el:fill, gr:pseudo} (we will think of our Stein manifolds as having 
boundaries), since the complex tangencies to the boundary induce a tight contact 
structure on the 3-manifold. For our purposes the exact definition of Stein 
manifold is unimportant (the curious reader is referred to 
\cite{el:stein}) as the following theorem gives a useful 
characterization of these manifolds.  But first we recall that a 
\dfn{Legendrian} knot in a contact 3-manifold $(M,\xi)$ is a curve
$\gamma:S^{1}\to M$ with $\gamma'(t)$ in $\xi_{t}$ for all $t\in 
S^{1}.$  To a Legendrian knot $\gamma$ bounding a surface $\Sigma$ we 
can assign two invariants.  The \dfn{Thurston-Bennequin invariant} of 
$\gamma,$ $\tb(\gamma,\Sigma),$ is simply the integer given by the 
framing induced on $\gamma$
by $\xi$ (where $\Sigma$ defines the zero framing).  The \dfn{rotation 
number} of $\gamma,$ $r(\gamma),$ is defined as follows: pick a trivialization
of $\xi\vert_\Sigma$ and let $T$ be a vector field tangent to $\gamma$.  Define
$r(\gamma,\Sigma)$ to be the degree of $T$ with respect to this trivialization.
The following theorem is implicit in Eliashberg's paper \cite{el:stein}.  For a complete
discussion of this theorem see the paper \cite{g:stein} of Gompf.

\begin{thm}\label{4Dstein}
	An oriented 4-manifold $X$ is a Stein manifold if and only if it has a handle
	decomposition with all handles of index less than or equal to 2 and 
	each 2-handle is attached to a Legendrian circle $\gamma$ with the framing 
	on $\gamma$ equal to $\tb(\gamma)-1$.
	Moreover, the first Chern class $c_1(X)$ is represented by the cocycle
	$$c=\sum r(\gamma_i)f_{h_i},$$
	where the sum is over the knots $\gamma_i$ to which the 2-handles $h_i$ are 
	attached
	and $f_{h_i}$ is the cochain that is 1 on core of $h_i$ and 0 elsewhere.
\end{thm}

\begin{figure}[ht]
		{\centerline{\epsfbox{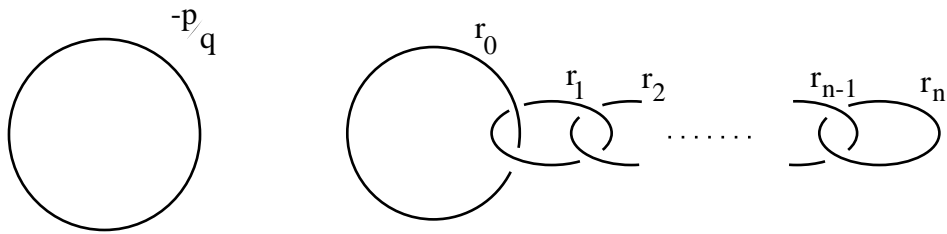}}}
		\caption{Two Kirby diagrams of  $L(p,q).$}
		\label{fig:lensexamples}
\end{figure}

One can use this to construct tight contact structures on every lens 
space $L(p,q)$ and compute their Euler class (since 
$e(\xi)=c_{1}(X)\vert_{\partial X}$).  To do this assume $p$ and $q$ are both positive 
(this is no restriction) and let $r_{0},r_{1},\ldots,r_{n}$ be a 
continued fraction expansion of $-\frac{p}{q}.$  The Kirby diagrams in 
Figure~\ref{fig:lensexamples} 
both represent $L(p,q)$ (for more on Kirby diagrams see 
\cite{gs}) and the diagram on the right can easily be made into an 
appropriate Legendrian link since all the surgery coefficients are 
integers less that 1.  In general we can only use this construction 
to construct one tight contact structure on $L(p,q),$ but on $L(p,1)$ 
we can do much better \cite{el:knots}.  If $p$ is odd, then we can 
realize all non-zero elements of $H^{2}(L(p,1);\Z)$ as Euler classes 
of a tight contact structure.  If $p$ is even, then we can realize 
all non-zero elements of 
$H^{2}(L(p,1);\Z)$ as ``$\Gamma(\xi)(s)$'' (where $s$ is 
as in Theorem~\ref{l-formula}) of a tight contact structure.  
We will see below that the ``missing 
classes'' above actually {\em cannot} be the Euler class 
(half Euler class) 
of a tight contact structure.  

In \cite{eg:tight} it was shown that 
one can do any Dehn surgery on unknots in $S^{3}$ (and some restricted surgeries 
on other knots) to obtain tight contact structures.  
Our main concern here is the existance of a tight contact structure 
realizing a particular half-Euler class.
\begin{thm}[\cite{eg:tight}]
Let $h$ be the 2-homology class determined by $h\cdot D=
\frac{1}{2} (q+1)\mod p.$
Any lens space $L(p,q)$ admits a tight contact structure 
with $e(\xi)=2h$ if $p$ is odd and $\Gamma(\xi)(s)=h$ if $p$ is even.
\end{thm}

\subsection{Uniqueness and Non-Existence of Tight Contact 
Structures}\label{main-section}

We begin by considering when the half-Euler class of a tight contact 
structure determines the structure.

\begin{thm}\label{lens}
	Let $L(p,q), p>0,$ be a lens space and $\xi_i$, $i=0,1$, be two tight contact 
	structures on $L(p,q)$.
	If 
	\begin{equation}
	\Gamma(\xi_i)(s')\cdot D=\pm \frac{1}{2}(q+1+p\Delta(s,s'))\mod p,
	\end{equation}
	for $i=0,1$,
	where $s, s'$ and $\Delta$ are as in \tm{l-formula},
	then $\xi_0$ and $\xi_1$ are contactomorphic.
\end{thm}

\bpr
    When $\Gamma(\xi_i)(s')\cdot D= \frac{1}{2}(q+1+p\Delta(s,s'))\mod p$ we can use 
    Formula~\ref{prac-l} to see that $l=-1+2np.$  Thus
    \Lm{simplify} allows us to arrange for $D_{\xi_{i}}$ to have exactly one 
    singular point which will have to be a positive elliptic point. 
    We can actually think of these isotopies as ambient isotopies of 
    $L(p,q).$  Thus we have isotoped the identity map to one which 
    takes a generalized projective plane (in the domain lens space) 
    with simple characteristic 
    foliation to a generalized projective plane (in the range lens 
    space) with simple characteristic foliation.
    We can further isotope our map so that it preserves the 
    characteristic foliation on $D.$  Thus \pr{foliation-determines} 
    will produce the desired contactomorphism.
    
    Every lens space $L(p,q)$ has an orientation preserving 
    diffeomorphism that acts on $H^{2}(L(p,q))$ by multiplication 
    by $-1$ \cite{bo}.  This allows us to reduce the
    $\Gamma(\xi_i,s')\cdot D= -\frac{1}{2}(q+1+p\Delta(s,s'))\mod p$ 
    case to the one above.
\epr

\bbr
  Some lens spaces have other orientation preserving 
  diffeomorphisms.  Using these one might hope to find other
  homology classes supporting at most one tight contact structure;  
  however, this does not seem to work since these
  diffeomorphisms permute the spin structures on $L(p,q).$  So the action of the 
  diffeomorphism on cohomology coupled with the action on the spin 
  structures conspire to prevent us from generalizing the above 
  theorem.
\eer

This theorem simplifies when $p$ is odd.

\begin{cor}
   Let $L(p,q), p>0,$ be a lens space and $\xi_i$, $i=0,1$, be two tight contact 
	structures on $L(p,q)$. If 
	$p$ is odd and 
	\begin{equation}
	e(\xi_i)(D)=\pm (q+1)\mod p,
	\end{equation}
	for $i=0,1$, then $\xi_0$ and $\xi_1$ are 
	contactomorphic. 
\end{cor}

We can also prove that some cohomology classes cannot be realized by 
any tight contact structure.

\begin{thm}
   	Let $\xi$ be a contact 
	structure on $L(p,q), p>1.$
	If 
	\begin{equation}
	\Gamma(\xi)(s')\cdot D=\pm \frac{1}{2}(q-1+p\Delta(s,s'))\mod p,
	\end{equation}
	where $s,s'$ and $\Delta$ are as in \tm{l-formula}, 
	then $\xi$ is overtwisted.
\end{thm}

\begin{figure}[ht]
		{\centerline{\epsfbox{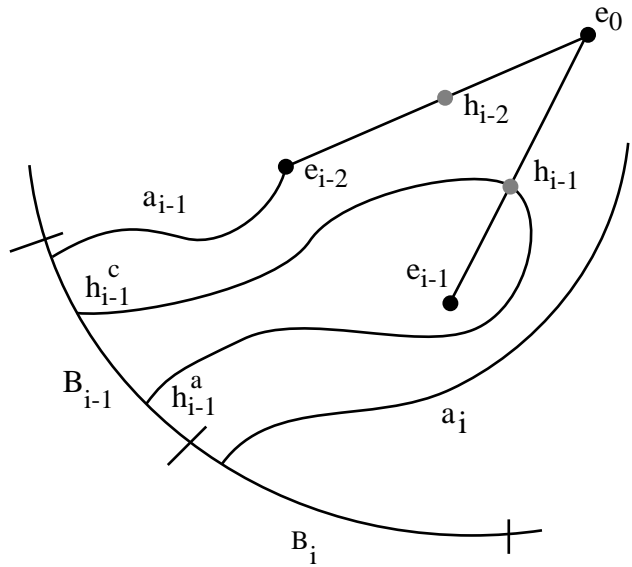}}}
		\caption{Impossible Configuration.}
		\label{bad}
	\end{figure}

\bpr
    We begin by assuming that $\xi$ is tight and proceed to find an 
    overtwisted disk.
    In the $\Gamma(\xi,s')\cdot D=\frac{1}{2}(q-1+p\Delta(s,s'))\mod p$ case we have 
    have $p$ elliptic points according to \Lm{l-formula} and 
    Remark~\ref{eh-formula}. Thus we can assume the graph of 
    singularities of $D_{\xi}$ is a star with $p-1$ branches. Using 
    the notation of Remark~\ref{setup} we can choose a point $x$ on 
    $C$ so that $h_{1}^{c}$ is closest (anti-clockwise) to $x$ and lying in 
    $B_{1}.$  The unstable seperatrix containing $h_{1}^{c}$ re-enters 
    $D$ after passing through $C$ in $p-1$ arcs, which we denote 
    $a_{2},\ldots, a_{p}.$  We now claim that $a_{i}$ ends at 
    $e_{i-1}$ for $i=2,\ldots (p-1).$  To see this assume it is false 
    and let $i$ be the 
    smallest index for which $a_{i}$ does not end at $e_{i-1}.$  
    If $a_{i}$ ends anti-clockwise of $e_{i-1}$ then both $h_{i-1}^{a}$ and 
    $h_{i-1}^{c}$ must be on $B_{i-1}$ contradicting 
    Remark~\ref{big-point}.  See Figure~\ref{bad}.
	Now if $a_{i}$ is clockwise of 
    $e_{i-1}$ then $h_{i-1}^{a},h_{i-1}^{c},\ldots,h_{p}^{a},h_{p}^{c}$ 
    must lie in $B_{i}\cup\ldots\cup B_{p}\cup B_{1}$  with no two 
    points coming from the same hyperbolic point lying in the same 
    $B_{j}.$  There are not, however, enough $B_{j}'s$ for this; thus 
    proving out claim.
    
\begin{figure}[ht]
		{\centerline{\epsfbox{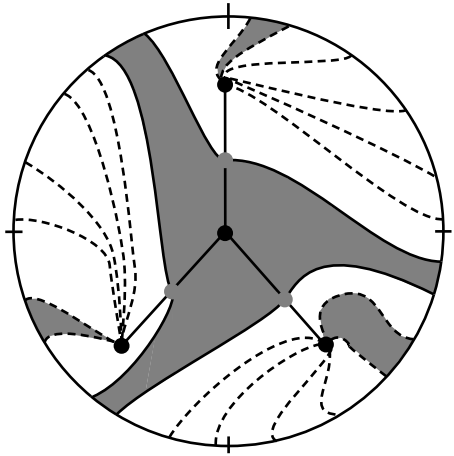}}}
		\caption{An Overtwisted Disk: dotted lines indicate arcs in the
			characteristic foliation corresponding to the continuation
			of unstable manifolds across $C,$ e.g. the $a_i$'s.}
		\label{fig:otd}
	\end{figure}

    Arguing in a similar fashion with the rest of the separatrices we 
    eventually see that the characteristic foliation must look like 
    the one shown in Figure~\ref{fig:otd}.  
    In this picture we can 
    explicitly find an overtwisted disk by canceling the 
    singularities on the boundary of the shaded region in 
    Figure~\ref{fig:otd}.
    
    Finally, if   $\Gamma(\xi,s')\cdot 
    D=-\frac{1}{2}(q-1+p\Delta(s,s'))\mod p$ then 
    one uses the diffeomorphism discussed at the end of the proof of 
    \tm{lens} to reduce to the above case.
\epr

Again we have a simpler statement when $p$ is odd.

\begin{cor}
   Let $\xi$ be a contact 
	structures on $L(p,q), p>1.$ If 
	$p$ is odd and 
	\begin{equation}
	  e(\xi)(D)=\pm (q-1)\mod p,
	\end{equation} 
	then $\xi$ is overtwisted.
\end{cor}

The above theorems do not provide a complete 
classification of contact structures on all lens spaces.  The best 
general statement that can be made is given in the following theorem.

\begin{thm}\label{short}
    On any lens space $L(p,q)$ there is at least one class in 
    $H^{2}(L(p,q))$ realized by a unique tight contact structure and 
    at least one class that cannot be realized by a tight contact structure.
\end{thm}

Currently,
contact structures are classified on
$L(p,q)$ only when  $p<4.$

\begin{thm}\label{complete}
  Classified up to isotopy
  \begin{enumerate}
    \item If $p=0$ then $L(p,q)=S^{1}\times S^{2}$ and there is a 
    unique tight contact structure.
    \item If $p=1$ then $L(p,q)=S^{3}$ and there is a unique tight 
    contact structure.
    \item If $p=2$ then $L(p,q)=\R P^{3}$ and there is a unique tight 
    contact structure.
    \item On $L(3,1)$ there are exactly two tight contact structures 
    (one for each non zero element in $H^{2}(L(3,1);\Z)$).
    \item On $L(3,2)$ there is exactly one tight contact structure 
    (realizing the zero class in $H^{2}(L(3,1))$).
  \end{enumerate}
\end{thm}

The first three statements were proved by Eliashberg \cite{el:twenty}.  
All but the first  
statement follow immediately from the theorems in this section.  
It is an interesting exercise to directly prove 4. by just considering 
$D_{\xi}$ and not resorting to the diffeomorphism used in \tm{lens}.  
Note that for all the examples mentioned in this theorem there is a 
unique contact structure up to contactomorphism but this is not always
the case, as exemplified by $L(4,1)$ which has at least two tight contact 
structures up to contactomorphism.

\subsection{Finiteness Results}

Though work of Kronheimer and Mrowka \cite{km} indicates that tight contact 
structures exist in only finitely many homotopy classes of plain 
fields, it is not, in general, known if any given 3-manifold has a finite 
number of tight contact structures.  There are examples of manifolds 
with infinitely many structures.  For example Giroux 
\cite{gi} and Kanda \cite{kanda} have shown that  $T^{3}$ has infinitely 
many tight contact structures.  To show there were infinitely many 
structures on $T^{3},$ essential use was made of incompressible 
tori in $T^{3}.$  One might hope that on atoroidal manifolds there are 
only finitely many tight contact structures. Thus lens spaces, being 
atoroidal, should have only a finite number of tight contact 
structures.  This is indeed the case.

\begin{thm}\label{finite}
  Any lens space admits only finitely many tight contact structures.
\end{thm}

\begin{proof}
  On $L(p,q)$ there are between $1$ and $p-1$ positive elliptic 
  singularities.  Once the number of positive elliptic singularities 
  is determined the entire characteristic foliation $D_{\xi}$ is determined by 
  the cyclic ordering of the $h_{i}^{c}$'s and $h_{i}^{a}$'s along 
  $C$ and the grouping of these points in the $B_{i}$'s.  Since there 
  are only a finite number of ways to order and group these points 
  the proof is complete.
\end{proof}

\bbr
  Refining the analysis in this proof of the structure of the 
  characteristic 
  foliation on $D$ one could derive a crude upper bound on the number of 
  tight contact structures on a given lens spaces.
\eer

{\em Acknowledgments:}
The author would like to thank Robert Gompf and \v Zarko Bi\v zaca for many helpful 
conversations concerning this paper, Robert Ghrist and Margaret 
Symington for providing many 
helpful comments on an early draft of this paper, and the University of Texas 
for their support of this work. The author also gratefully acknowledges support by an NSF 
Post-Doctoral Fellowship(DMS-9705949) and Stanford University during the writing of this paper.


\end{document}